%

\magnification=1200
\input amstex
\documentstyle{amsppt}

\define\s{{\frak s}}
\define\g{{\frak g}}
\define\bb{{\frak b}}
\define\dd{{\frak d}}
\define\up{{\Cal U}\text{-prod}\,\omega}
\define\upp{{\Cal U}'\text{-prod}\,\omega}
\define\cf#1{{\text{cf}\,(#1)}}
\define\cu{{\Cal U}}
\define\ai{\forall^\infty}
\define\ei{\exists^\infty}
\define\ww{{}^\omega\omega}
\topmatter
\title
On the Cofinality of Ultrapowers
\endtitle
\author
Andreas Blass\\and\\
Heike Mildenberger
\endauthor
\address
Mathematics Dept., University of Michigan, Ann Arbor, 
MI 48109, U.S.A.
\endaddress
\email
ablass\@umich.edu, mildenbe\@umich.edu
\endemail
\thanks 
The first author is partially supported by National Science Foundation
grant DMS-9505118 and the second by a postdoctoral fellowship from the
Deutsche Forschungsgemeinschaft.
\endthanks
\abstract
We prove some restrictions on the possible cofinalities of ultrapowers
of the natural numbers with respect to ultrafilters on the natural
numbers.  The restrictions involve two cardinal characteristics of the
continuum, the splitting number $\s$ and the groupwise density number
$\g$. 
\endabstract
\endtopmatter
\document

\head
1. Introduction
\endhead

All ultrafilters considered in this paper are non-principal
ultrafilters on the set $\omega$ of natural numbers.  We shall be
concerned with the possible cofinalities $\cf\up$ of ultrapowers of
$\omega$ with respect to such ultrafilters.  We shall show that no
cardinal below the groupwise density number $\g$ (see definition
below) can occur as such a cofinality and that at most one cardinal
below the splitting number $\s$ can so occur.  The proof for $\s$,
when combined with a result of Nyikos, gives the additional
information that all $P_{\bb^+}$-point ultrafilters are nearly
coherent.

In Section 2, we review the necessary terminology and some previously
known results.  In Section~3, we prove the result concerning $\g$.
Finally, in Section~4, we prove the result concerning $\s$, we show
that in the statement of that result ``at most one cardinal'' cannot
be improved to ``no cardinal,'' and we deduce the result about
$P_{\bb^+}$-points.  

We thank Simon Thomas for posing the question whether $\cf\up$ can
ever be smaller than $\g$.

\head
2. Preliminaries
\endhead

We write $\ei$ and $\ai$ for the quantifiers ``there exist infinitely
many'' and ``for all but finitely many,'' respectively.  Any
ultrafilter (by which we always mean a non-principal ultrafilter on
$\omega$) $\cu$ will also be used as a quantifier meaning ``for almost
all with respect to $\cu$,'' i.e.,
$$
(\cu n)\,\varphi(n)\iff\{n\mid\varphi (n)\}\in\cu.
$$
Thus, the quantifier $\cu$ is intermediate between $\ai$ and $\ei$ in
the sense that $(\ai n)\,\varphi(n)\implies(\cu n)\,\varphi(n)\implies
(\ei n)\,\varphi(n)$ for any predicate $\varphi$ on natural numbers.

The ultrapower $\up$ is formed from the set $\ww$ of all functions
$f:\omega\to\omega$ by identifying $f$ with $g$ whenever $(\cu
n)\,f(n)=g(n)$.  It is linearly ordered by the relation
$$
f\leq_\cu g\iff(\cu n)\,f(n)\leq g(n).
$$
By $\cf\up$ we mean the cofinality of this ordering, the smallest
cardinality of a subset $\Cal C$ of $\ww$ such that every $f\in\ww$ is
$\leq_\cu$ some $g\in\Cal C$.

This cofinality obviously satisfies $\bb\leq\cf\up\leq\dd$, where the
bounding number $\bb$ and the dominating number $\dd$ are defined as
follows.  (For more information on these and other cardinal
characteristics of the continuum, see the survey papers \cite{7,11}.)
$\dd$ is the minimum size of a family $\Cal D\subseteq\ww$ such that,
for each $f\in\ww$ there is some $g\in\Cal D$ satisfying $(\ai
n)\,f(n)\leq g(n)$.  The definition of $\bb$ is the same except that
$\ai$ is replaced with $\ei$.

In addition to $\bb$ and $\dd$, three other cardinal characteristics
of the continuum, $\s$, $\g$, and cov(B), will play a role in this
paper.  

The splitting number $\s$ is defined as the minimum size of a
family $\Cal S$ of subsets of $\omega$ such that every infinite
$X\subseteq\omega$ is split by some $Y\in\Cal S$ in the sense that
both $X\cap Y$ and $X-Y$ are infinite.  

To define $\g$, we first need the notion of groupwise density.  A
family $\Cal G$ of infinite subsets of $\omega$ is said to be
groupwise dense if it is closed under infinite subsets and finite
modifications and if, whenever $\omega$ is partitioned into finite
intervals, the union of some infinitely many of these intervals is in
$\Cal G$.  Then $\g$ is defined as the minimum number of groupwise
dense families with empty intersection.  (See \cite3 for more
information about groupwise density and $\g$.)

Finally, cov(B) is defined to be the minimum number of meager sets
(i.e., sets of the first Baire category) needed to cover the real
line. 

We shall be concerned with restrictions, in terms of cardinal
characteristics of the continuum, on the possible values of $\cf\up$.
The following theorem of Canjar \cite{4,5} and Roitman \cite9 suggests
that the trivial restriction $\bb\leq\cf\up\leq\dd$ is all one can
hope for.

\proclaim{Theorem 1 \cite{4,5,9}}
It is consistent (relative to ZFC) that $\bb\ll\dd$ and every regular
cardinal $\kappa$ in the range $\bb\leq\kappa\leq\dd$ occurs as
$\cf\up$ for some $\cu$.
\endproclaim

The model used to prove this theorem is the Cohen model, obtained by
adding a large number of Cohen-generic reals to any model of ZFC.  We
shall see that the trivial lower bound $\bb$ for all $\cf\up$ can be
improved in some models (but not in all, by Theorem~1).

Canjar also showed that the trivial upper bound $\dd$ cannot be
improved in any model where $\dd$ is regular.  

\proclaim{Theorem 2 \cite6}
There exists an ultrafilter $\cu$ with $\cf\up=\cf\dd$.  In
particular, if $\dd$ is regular then it occurs as $\cf\up$ for some
$\cu$. 
\endproclaim

For any ultrafilter $\cu$ and any function $f:\omega\to\omega$, the
image $f(\cu)$ is defined as the ultrafilter $\{X\subseteq\omega\mid
f^{-1}(X)\in\cu\}$.  (Contrary to our convention, this may be
a principal ultrafilter, but only if $f$ is constant on some set in
$\cu$; we shall use  $f(\cu)$ only for finite-to-one functions $f$, so
no real difficulty arises.)  Two ultrafilters $\cu$ and $\cu'$ are
said to be nearly coherent if $f(\cu)=f'(\cu')$ for some finite-to-one
functions $f$ and $f'$.  It is shown in \cite1 that the same relation
of near coherence would be obtained if we required in the definition
that $f=f'$ and that $f$ be monotone.  It is also shown there that
near coherence is an equivalence relation and that, whenever $\cu$ and
$\cu'$ are nearly coherent, then $\cf\up=\cf\upp$ (because both of
these ultrapowers have cofinal submodels isomorphic to
$f(\cu)\text{-prod}\,\omega$).  The principle of near coherence of
filters (NCF), introduced in \cite1 and proved consistent in \cite2,
asserts that every two non-principal ultrafilters on $\omega$ are
nearly coherent.

\head
3. Groupwise Density Gives a Lower Bound
\endhead

In this section, we prove the following answer to a question raised by
Simon Thomas (private communication).

\proclaim{Theorem 3}
For every non-principal ultrafilter $\cu$ on $\omega$,
$\cf\up\geq\g$. 
\endproclaim

\demo{Proof}
Suppose $\Cal C\subseteq\ww$ is cofinal with respect to $\leq_\cu$.
We shall associate to each $f\in\Cal C$ a groupwise dense family $\Cal
G_f$ in such a way that the intersection of these families is empty.
Thus, we shall have $\g\leq|\Cal C|$, which establishes the theorem.

By increasing them if necessary, we may assume without loss of
generality that all the functions $f\in\Cal C$ satisfy $f(n)\geq n$
for all $n$.  To define $\Cal G_f$, we first define, for each infinite
$X\subseteq\omega$, the function $\nu_X:\omega\to\omega$ sending each
natural number $n$ to the next larger element of $X$.  Then let 
$$
\Cal G_f=\{X\subseteq\omega\mid X\text{ is infinite and
}f<_\cu\nu_X\}
$$
for each $f\in\Cal C$.  Since these $f$'s are cofinal in $\up$, the
intersection of the corresponding $\Cal G_f$'s must be empty.  It is
also clear that each $\Cal G_f$ is closed under infinite subsets and
under finite modifications.  So to verify that each $\Cal G_f$ is
groupwise dense, thus completing the proof, it remains only to check
that, if $f$ is fixed and if $\omega$ is partitioned into finite
intervals then the union of some infinitely many of these intervals is
in $\Cal G_f$.

Inductively select intervals $I_k$ from the given partition so that
the first element of $I_{k+1}$ is greater than $f(x)$ for all $x\in
I_k$ and all smaller $x$.  Let $X$ be the union of the even-numbered
intervals, $I_{2j}$, and $Y$ the union of the odd-numbered ones.  

For any natural number $p$ in the interval $(\max I_{n-1}, \max I_n]$,
one of $\nu_X(p)$ and $\nu_Y(p)$ (depending on the parity of $n$) will
be $\min I_{n+1}$, which is greater than $f(p)$.  Thus, every natural
number $p$, except for the finitely many below $\max I_0$, is in one
of the two sets $\{n\in\omega\mid f(n)<\nu_X(n)\}$ and
$\{n\in\omega\mid f(n)<\nu_Y(n)\}$.  Therefore, one of these sets is
in $\cu$, which means that one of $X$ and $Y$ is in $\Cal G_f$.  Since
both $X$ and $Y$ are unions of infinitely many intervals from the
given partition, this completes the proof that $\Cal G_f$ is
groupwise dense and thus completes the proof of the theorem.
\qed\enddemo

It is well-known (see \cite3) that $\g\leq\dd$.  The following
corollary gives an improvement when $\dd$ is singular.

\proclaim{Corollary 4}
$\g\leq\cf\dd$.
\endproclaim

\demo{Proof}
Combine Theorems~2 and 3.
\qed\enddemo

Encouraged by Theorem~3, one might look for additional cardinal
characteristics that give lower bounds on the possible cofinalities of
$\up$.  Such characteristics must be $\leq\dd$ and, to avoid
trivialities, $\not\leq\bb$.  Inspection of the diagrams of cardinal
characteristics in \cite{11} provides just two such characteristics, the
splitting number $\s$ and the covering number for category cov(B).
(If one counts the somewhat 
artificial $\min\{\frak r,\dd\}$ as a characteristic, then it also
lies in the desired region.  The following remark about cov(B) applies
to it as well.)  If we add a large number $\kappa$ of Cohen reals to a
model of set theory, then the resulting model has cov(B) large but
has, by the proof of Theorem~1, ultrafilters with $\cf\up=\aleph_1$.
So cov(B) cannot serve as a lower bound for $\cf\up$.  That leaves
$\s$ as a possibility, which we analyze in the next section.

\head
4. The Splitting Number
\endhead

Unlike $\g$, the splitting number $\s$ is not in general a lower bound
for $\cf\up$.  The proof involves the notion of (pseudo-)$P_\kappa$
point.  An ultrafilter $\cu$ is called a $P_\kappa$ point if, for
every family $\Cal F\subseteq\cu$ with $|\Cal F|<\kappa$, there is some
$A\in\cu$ with $A-F$ finite for all $F\in\Cal F$.  Pseudo-$P_\kappa$
points are defined similarly, except that $A$ is not required to be in
$\cu$, only to be infinite.  We shall need the following results of
Nyikos, folklore, and Shelah, respectively.  (Although Nyikos's paper
\cite8 is not yet published, Proposition~5 and its proof were in a
1984 letter from Nyikos to the first author.)

\proclaim{Proposition 5 \cite8}
If $\cu$ is a pseudo-$P_\kappa$ point and $\kappa>\bb$, then
$\cf\up=\bb$.
\endproclaim

\proclaim{Proposition 6}
If $\cu$ is a pseudo-$P_\kappa$ point then $\s\geq\kappa$.
\endproclaim

\proclaim{Proposition 7 \cite2}
It is consistent relative to ZFC that $\bb=\aleph_1$ and there is a
$P_{\aleph_2}$-point.
\endproclaim

Since the first two of these propositions are fairly easy, we give
their proofs.  For Proposition~7, we refer to Theorem~6.1 of \cite2,
which gives (more than) a model with a $P_{\aleph_2}$-point and
another ultrafilter generated by $\aleph_1$ sets.  The latter gives us
$\bb=\aleph_1$ because, by a theorem of Solomon \cite{10}, no
ultrafilter can be generated by fewer than $\bb$ sets.

\demo{Proof of Proposition 5}
Let $\cu$ be a pseudo-$P_\kappa$ point with $\kappa>\bb$, and let
$\Cal C\subseteq\ww$ be a family of cardinality $\bb$ such that for every
$f\in\ww$ there is $g\in\Cal C$ with $(\ei n)\,f(n)\leq g(n)$.  By
increasing each $g\in\Cal C$ if necessary, we can assume that $g$ is a
monotone non-decreasing function.  To complete the proof, we show that
$\Cal C$ is cofinal with respect to the linear ordering $\leq_\cu$ of
$\up$.  

Suppose to the contrary that $h\in\ww$ is such that $g\leq_\cu h$ for
all $g\in\Cal C$.  This means that the sets $M_g=\{n\in\omega\mid
g(n)\leq h(n)\}$ are in $\cu$ for all $g\in\Cal C$.  Since $|\Cal
C|=\bb<\kappa$ and since $\cu$ is a pseudo-$P_\kappa$ point, there is
an infinite set $X\subset\omega$ such that each $X-M_g$ is finite.  As
in the proof of Theorem~3, let $\nu_X(n)$ denote the next member of
$X$ after $n$.  By our original choice of $\Cal C$, there is $g\in\Cal
C$ such that $h(\nu_X(n))<g(n)$ for infinitely many $n$.  For each
such $n$ we have, since $g$ is non-decreasing,
$h(\nu_X(n))<g(\nu_X(n))$ and therefore $\nu_X(n)\in X-M_g$.  But this
applies to infinitely many $n$, giving infinitely many $\nu_X(n)$,
contrary to the fact that $X-M_g$ is finite.
\qed\enddemo

\demo{Proof of Proposition 6}
Let $\cu$ be a pseudo-$P_\kappa$ point and let $\Cal S$ be a family of
fewer than $\kappa$ subsets of $\omega$.  We must find an infinite set
$X\subseteq\omega$ that is not split by any member of $\Cal S$.

For each $Y\in\Cal S$, let $Y'$ be $Y$ or $\omega-Y$, whichever is in
$\cu$.  As $\cu$ is a pseudo-$P_\kappa$ point, there is an infinite
$X$ such that $X-Y'$ is finite for all $Y\in\Cal S$.  This $X$ is
clearly not split by any such $Y$.
\qed\enddemo

\proclaim{Corollary 8}
It is consistent, relative to ZFC, that there is a non-principal
ultrafilter $\cu$ on $\omega$ with $\cf\up<\s$.
\endproclaim

\demo{Proof}
In the model given by Proposition~7, let $\cu$ be a $P_{\aleph_2}$
point.  Its existence gives $\s\geq\aleph_2$ by Proposition~6, and we
also have, by Propositions~5 and 7, $\cf\up=\bb=\aleph_1$.
\qed\enddemo

Although Corollary 8 shows that it is consistent for the set of
cofinalities of ultrapowers of $\omega$ to contain a cardinal below
$\s$, we shall see that this set cannot contain two cardinals below
$\s$.  That will be a consequence of the following theorem.

\proclaim{Theorem 9}
Suppose $\cu$ and $\cu'$ are non-principal ultrafilters on $\omega$
such that both $\cf\up$ and $\cf\upp$ are smaller than $\s$.  Then
$\cu$ and $\cu'$ are nearly coherent.
\endproclaim

\demo{Proof}
Let $\cu$ and $\cu'$ satisfy the hypotheses of the theorem, and
suppose these ultrafilters are not nearly coherent.  Let $\Cal C$ and
$\Cal C'$ be subfamilies of $\ww$, each of size $<\s$, and cofinal
with respect to $\leq_\cu$ and $\leq_{\cu'}$ respectively.  Let $\Cal
D$ be the set of functions of the form $\max\{g,g'\}$, where $g\in\Cal
C$, $g'\in\Cal C'$, and max means the pointwise maximum of the
functions.  Then, for each $f\in\ww$, there is an $h\in\Cal D$ such
that both inequalities $f\leq_\cu h$ and $f\leq_{\cu'}h$ hold.

Temporarily fix some $h\in\Cal D$.  Partition $\omega$ into finite
intervals $I_n=[a_n,a_{n+1})$ such that $h(x)<a_{n+1}$ for all
$x<a_n$.  (It is trivial to produce such $a_0=0<a_1<a_2<\dots$
inductively.)  Let $p:\omega\to\omega$ be the function that sends all
points in $I_n$ to $n$, for all $n$.  Since $p$ is finite-to-one and
since $\cu$ and $\cu'$ are not nearly coherent, the ultrafilters
$p(\cu)$ and $p(\cu')$ are distinct, so one contains a set whose
complement is in the other.  Pulling these sets back along $p$, we get
two sets, say $A\in\cu$ and $A'\in\cu'$, each a union of some $I_n$'s,
but with no $I_n$ in common.

Define $q(x)=p(x)+1$.  Applying again the fact that $\cu$ and $\cu'$
are not nearly coherent, we have $q(\cu)\neq p(\cu')$, so we can get a
set in $q(\cu)$ whose complement is in $p(\cu')$.  Pulling these sets
back along $q$ and $p$ respectively, we get $B\in\cu$ and $B'\in\cu'$,
each a union of some $I_n$'s, and such that we never have an
$I_n\subseteq B$ and $I_{n+1}\subseteq B'$.

Arguing analogously with $p(\cu)\neq q(\cu')$, we get $C\in\cu$ and
$C'\in\cu'$, each a union of some $I_n$'s, such that we never have an
$I_n\subseteq C'$ and $I_{n+1}\subseteq C$.

Let $D=A\cap B\cap C$ and $D'=A'\cap B'\cap C'$.  Then $D\in\cu$,
$D'\in\cu'$, and both are unions of some $I_n$'s.  Furthermore, if a
particular $I_n$ is included in $D$ then neither it nor its neighbors
$I_{n\pm1}$ can be included in $D'$.

Let $E$ be the union of all the $I_n$'s and $I_{n+1}$'s such that
$I_n\subseteq D$, i.e., the union of the intervals that constitute $D$
and their right neighbor intervals.  Define $E'$ similarly from $D'$,
and note that $E$ and $E'$ are disjoint.

I claim that, if $X$ is an infinite subset of $\omega$ and if
$\nu_X\leq_\cu h$, then $X\cap E$ is infinite.  To see this, notice
first that the set $\{k\in\omega\mid\nu_X(k)\leq h(k)\}$, being in
$\cu$, must contain infinitely many points $k\in D$ because
$D\in\cu$.  For each of these infinitely many $k$, there is an element
of $X$, namely $\nu_X(k)$, in the interval $[k,h(k)]$.  By our choice
of the intervals $I_n$, this element of $X$ is either in the same
interval as $k$ or in its right neighbor.  In either case, it is in
$E$ because $k\in D$.  Thus, we have infinitely many (since $k$ can be
arbitrarily large) elements of $X\cap E$, as claimed.

Similarly, if $\nu_X\leq_{\cu'}h$, then $X\cap E'$ is infinite and
therefore so is $X-E$ since $E$ and $E'$ are disjoint.

Now un-fix $h$.  For each $h\in\Cal D$, the preceding discussion
produces an $E$, which we now call $E_h$ to indicate its dependence on
the (previously fixed) $h$.  For any infinite subset $X$ of $\omega$,
the function $\nu_X$ is majorized, with respect to both $\leq_\cu$ and
$\leq_{\cu'}$, by some $h\in\Cal D$.  then the preceding discussion
shows that $X$ is split by the corresponding $E_h$.  Therefore,
$\{E_h\mid h\in\Cal D\}$ is a splitting family.  But this is absurd,
as $|\Cal D|<\s$.
\qed\enddemo

\proclaim{Corollary 10}
At most one cardinal smaller than $\s$ can occur as $\cf\up$.
\endproclaim

\demo{Proof}
Combine Theorem 9 and the fact that nearly coherent ultrafilters
produce ultrapowers of the same cofinality.
\qed\enddemo

\proclaim{Corollary 11}
Any two pseudo-$P_{\bb^+}$ points are nearly coherent.
\endproclaim

\demo{Proof}
If two ultrafilters are pseudo-$P_{\bb^+}$ points, then the
corresponding ultrapowers have cofinality $\bb$ by Proposition~5, and
this is smaller than $\s$ by Proposition~6.  So Theorem~9 applies and
gives the required near coherence.
\qed\enddemo

\remark{Remark}
For an ultrafilter $\cu$ to have a small system of generators and for
its ultrapower $\up$ to have small cofinality are in some sense
antithetical properties.  Specifically, the proof of Theorem~16 in
\cite1 shows that the number of generators of $\cu$ and $\cf\up$
cannot both be smaller than $\dd$.  Yet each property, when it holds
of two ultrafilters (with an appropriate sense of ``small'') implies
near coherence.  For $\cf\up$, the appropriate sense of ``small'' is
$<\s$ and the relevant result is Theorem~9 above.  For the number of
generators of $\cu$, the appropriate sense of ``small'' is $<\dd$, for
Corollary~13 of \cite1 says that any two ultrafilters generated by
fewer than $\dd$ sets are nearly coherent.
\endremark

\enddocument